\documentclass[12pt,leqno]{article}
\usepackage{amsmath,amsthm,amssymb}

\def\init{\setcounter{equation}{0}}
\newtheorem{theorem}{Theorem}[section]

\newcommand{\e}{{\varepsilon}}

\newcommand{\rw}{\rightarrow}

\usepackage{tikz}
\usetikzlibrary{decorations.pathreplacing}

\usepackage{verbatim}

\title{Behavior  of null-geodesics in the interior of Reissner-Nordstrom  black hole 
}

\author{G.Eskin, \ \ \  Department of Mathematics, UCLA,\\ Los Angeles,
CA 90095-1555, USA. \ E-mail: eskin@math.ucla.edu
}

\begin{document}

\maketitle

\begin{abstract}
We show  that an incoming  null-geodesic   belonging  to a plane  passing  through the origin  
and starting outside the outer  horizon  crosses the outer  and the inner horizons.  Then   it turns at some  point  inside  
the inner  horizon  and    approaches the inner horizon  when the time tends  to the  infinity.
We also construct  a geometric optics   solution  of the Reissner-Nordstrom  equation  that has support  in a neighborhood  of the 
null-geodesic.
\end{abstract}

{\bf Keywords.}  Reissner-Nordstrom black hole,  null-bicharacteristics,\\
\indent  geometric optics solutions.

\section{Introduction}
\init

The  Reissner-Norstrom metric  (cf. [7],  [9],  [12])  is a spherically symmetric  metric  having the following form  in Cartesian 
coordinates  (cf.  [11]):
\begin{equation}															\label{eq:1.1}
ds^2=\Big(1-\frac{2m}{r}+\frac{e^2}{r^2}\Big)dx_0^2-dx_1^2-dx_2^2-dx_3^2-2\Big( \frac{2m}{r}-\frac{e^2}{r^2}      \Big)dx_0 dr
-\Big(\frac{2m}{r}-\frac{e^2}{r^2}\Big)dr^2,
\end{equation}
where $r=\sqrt{x_1^2+x_2^2+x_3^2}$.

Let
\begin{equation}															\label{eq:1.2}
\Box_gu=0
\end{equation}
  be the wave  operator  corresponding  to  (\ref{eq:1.1}).

The symbol  (Hamiltonian)  of  $\Box_g$  has the form
\begin{equation}															\label{eq:1.3}
H=\xi_0^2-\sum_{j=1}^3\xi_j^2+\Big(\frac{2m}{r}-\frac{e^2}{r^2}\Big)\Big(-\xi_0+\sum_{j=1}^3x_j\xi_j\Big)^2,
\end{equation}
where  $x_0$   is the time coordinate,  $x=(x_1,x_2,x_3)$
are  the space  coordinates,  $(\xi_0,\xi_1,\xi_2,\xi_3)$  are dual  to  $(x_0,x_1,x_2,x_3)$,  $r=\sqrt{\sum_{j=1}^3 x_j^2}$.

Denote
\begin{equation}															\label{eq:1.4}
f=1-\frac{2m}{r}+\frac{e^2}{r^2}.
\end{equation}
Note that  $f=0$ has  two  real roots
\begin{equation}													   		\label{eq:1.5}
r_+=m+\sqrt{m^2-e^2},\ \ \ r_-=m-\sqrt{m^2-e^2},\ \ \ f=\frac{(r-r_+)(r-r_-)}{r^2},
\end{equation}
assuming that  $e^2<m^2$.

It  follows   from  (\ref{eq:1.1})   that  $r=r_+$  and  $r=r_-$  are   the outer  and inner  horizons  for  
the  Reissner-Nordstrom (RN) metric.

Sometimes in the case $e^2<m^2$  the RN  metric  is called  sub-extremal RN  metric.

If $e^2=m^2$  then  $r_+=r_-=m$  and  $f=\frac{(r-m)^2}{r^2}$.   In this case  RN  metric is called  extremal.  It has  only one    
event horizon   $r=m$.
In the case  $e^2>m^2$  we have  $f>0$  for all  $0<r<+\infty$  and    the  RN  metric  has  no  event horizon and is called  the naked 
singularity RN  metric.

The Reissman-Nordstrom (RN)  metric is a solution  of  
the Einstein  equations  of general relativity  in vacuum.  It was investigated in many papers
with the emphasis  on the study  of the stability  of the metric  among the solutions  of the Einstein  equations  (cf.  [1],  [5],  [8],  
[10]).

We will not use  that  (\ref{eq:1.1})  is a solution of the Einstein  equations.

The null-geodesics  of (\ref{eq:1.1})  are the projections on the $x$-space of  the  null-bicharacteristics (see \S 2  for the details).

We shall give   an explicit description of the behavior of  null-geodesics  of
the  wave equation  $\Box_gu=0$  and we shall construct  a geometric  optics  solution  of  (\ref{eq:1.2})  
with the support modulo  lower  order terms   in a neighborhood  of a null-geodesic.  We are  not  using  the expansion  of the solution  
in spherical  harmonics.  Instead  we are constructing  the solution   in three-dimensional  case using explicit formulas  for  
the  solutions  in    two-dimensional slices  (cf. [3]).  This way we  obtain  a more  explicit  formula than when using the spherical 
coordinates.  In particular,  we show that our  solution  has  a finite  limit  when  $x_0\rw +\infty$  along  any  null-geodesic
approaching  the inner horizon.

The plan of the paper is the following:

In \S 2  we study  the behavior   of null-geodesics  located on a plane passing through the origin.  We show  that  any incoming 
geodesic  starting  outside the black hole  crosses  the outer event  horizon $r=r_+$   when  the  time  $x_0$  is increases.  Then it
crosses also  the inner horizon  $r=r_-$  and reaches  a turning point $r=r_0<r_-$  where  $r_0$  depends  on the initial data   of
the null-bicharacteristic.   
After this it
 approaches  the inner horizon   from the inside
when  $x_0\rw +\infty$.

In \S3  we  construct  geometric  optics  solutions  $u_N$  of  (\ref{eq:1.2})  depending on  a large parameter  $k$.  This  geometric  optics  
solution  has  a support   in a neighborhood  of some null-geodesics $\gamma_0$.   Thus  the time evolution  of geometric  optics  solution   
follows  modulo lower order terms   the time evolution of  the  null-geodesic  $\gamma_0$.

In \S4  we study  the cases  of  extremal and naked  singularity  RN  metrics.

In \S5  we summarize  the results  of the paper.

\section{The behavior  of null-bicharacteristics  in the case  of Reissner-Nordstrom  metric}
\init

The  equations of  null-bicharacteristics  in Cartesian coordinates has the form 
\begin{align}															\label{eq:2.1}
&\frac{dx_k}{ds}=\frac{\partial H}{\partial \xi_k},\ \ x_k(0)=y_k,\ \ 0\leq k \leq 3,
\\
\nonumber
&\frac{d\xi_k}{ds}=-\frac{\partial H}{\partial x_k},\ \ \xi_k(0)=\eta_k,\ \ 0\leq k \leq 3,
\end{align}
where  $(\xi_0,\xi_1,\xi_2,\xi_3)$  are dual  coordinates  to  $(x_0,x_1,x_2,x_3)$,   $H$  is the same as in (\ref{eq:1.3}).  We study  
the restriction of the metric (\ref{eq:1.1}) to the plane  $x_3=0$.
Consider  the null-bicharacteristic such that  $x_3(0)=0,\xi_3(0)=0$.  We have
\begin{align}															\label{eq:2.2}
&\frac{dx_3}{ds}=\frac{\partial H}{\partial \xi_3} =-2\xi_3+    \Big(\frac{2m}{r}-\frac{e^2}{r^2}\Big)2\Big(-\xi_0+\sum_{k=1}^3x_k\xi_k\Big)x_3
\\
\nonumber
&\frac{d\xi_3}{ds}=-\frac{\partial H}{\partial x_3}  
 = - \Big(\frac{2m}{r}-\frac{e^2}{r^2}\Big)2\Big(-\xi_0+\sum_{k=1}^3x_k\xi_k\Big)\xi_3
\\
\nonumber
&+\Big(\frac{2m}{r^2}\frac{x_3}{r}+\frac{2e^2}{r^3}\frac{x_3}{r}\Big)\Big(-\xi_0+\sum_{k=1}^3x_k\xi_k\Big)^2.
\end{align}
Since  the initial data  
$x_3(0)=0,\xi_3(0)=0$ are zero,  by the uniqueness  theorem  for the  system  (\ref{eq:2.2}),   we have  that  
$x_3(s)\equiv 0,\xi_3(s)\equiv 0$,  i.e.
the null-bicharacteristic stays  in the plane  $x_3=0$   and   $\xi_3=0$.   Therefore  we can  restrict  the null-bicharacteristic  to the plane 
$x_3=0$  with  $\xi_3=0$  and  the restricted  Hamiltonian has the form
$$
H_0(x_1,x_2,\xi_0,\xi_1,\xi_2)=\xi_0^2-\sum_{k=1}^2\xi_k^2+\Big(\frac{2m}{\rho}-\frac{e^2}{\rho^2}\Big)\Big(-\xi_0+\sum_{k=1}^2x_k\xi_k\Big)^2,
$$
where
$\rho=\sqrt{x_1^2+x_2^2}$.   In polar  coordinates   $(\rho,\varphi)$  the Hamiltonian has the form
\begin{align}																\label{eq:2.3}
H_0(\rho,\varphi,\xi_0,\xi_\rho,\xi_\varphi)&=
\xi_0^2-\xi_\rho^2-\frac{1}{\rho^2}\xi_\varphi^2+\Big(\frac{2m}{\rho}-\frac{e^2}{\rho^2}\Big)(-\xi_0+\xi_\rho)^2
\\
\nonumber
&=(2-f)\xi_0^2+2(f-1)\xi_0\xi_\rho-f\xi_\rho^2-\frac{1}{\rho^2}\xi_\varphi^2,
\end{align}
where $f$  is the same as in (\ref{eq:1.4}).  The  system of  null-bicharacteristics has  the form
\begin{align}																\label{eq:2.4}
&\frac{d\rho}{ds}=\frac{\partial H_0}{\partial \xi_\rho}=2(f-1)\xi_0-2f\xi_\rho,\ \ \rho(0)=\rho_0,
\\
\nonumber
&\frac{d\varphi}{ds}=\frac{\partial H_0}{\partial\xi_\varphi}=-\frac{2\xi_\varphi}{\rho^2},\ \ \ \varphi(0)=\varphi_0,
\ \ \ \frac{d\xi_0}{ds}=-\frac{\partial H_0}{\partial x_0},\ \ \xi_0(0)=\xi_0^{(0)},
\\
\nonumber
&\frac{d\xi_\rho}{ds}=-\frac{\partial H_0}{\partial\rho},\ \ \ \xi_\rho(0)=\xi_\rho^{(0)},
\ \ \ \frac{d\xi_\varphi}{ds}=-\frac{\partial H_0}{\partial \varphi},\ \ \xi_\varphi(0)=\xi_\varphi^{(0)},
\\
\nonumber
&\frac{dx_0}{ds}=\frac{\partial H_0}{\partial \xi_0}=2(2-f)\xi_0+2(f-1)\xi_\rho,\ \ x_0(0)=x^{(0)}.
\end{align}

Since  $H_0$  is  independent  of $x_0$  and $\varphi$   we have   that   $\xi_0(s)=\xi_0^{(0)},\xi_\varphi(s)=\xi_\varphi^{(0)}$  for  all $s$.
For the simplicity of notations  we shall write  $\xi_0,\xi_\varphi$  instead of  $\xi_0^{(0)},\xi_\varphi^{(0)}$.
 Denote  by  $\gamma_0$  the null-bicharacteristic with the initial data
$\rho_0,\varphi_0,\xi_{\rho}^{(0)},\xi_{\varphi}^{(0)}, \xi_0^{(0)}$.
Note that  the null-bicharacteristic mean  that
\begin{equation}																\label{eq:2.5}
H_0(\rho(s),\varphi(s),\xi_0(s),\xi_\rho(s),\xi_\varphi(s))=0 \ \ \mbox{for all}\ \ s\geq 0.
\end{equation}
Therefore  we can find  $\xi_\rho$  from   (\ref{eq:2.5})  (cf. [3]).
We  get
\begin{multline}																	\label{eq:2.6}
\xi_\rho^\pm=\xi_\rho^\pm(\rho,\xi_0,\xi_\varphi)=
\frac{                   
-(f-1)\xi_0\pm\sqrt{(f-1)^2\xi_0^2+f\big[(2-f)\xi_0^2-\frac{1}{\rho^2}\xi_\varphi^2\big]}
}
{-f}
\\
=\frac{         
-(f-1)\xi_0\pm\sqrt{\xi_0^2-f\frac{\xi_\varphi^2}{\rho^2}}
}
{-f}
\end{multline}
In particular,  when $s=0$  we have
\begin{equation}																\label{eq:2.7}
\xi_\rho^\pm(0)=\xi_\rho^\pm(\rho_0,\xi_0,\xi_\varphi).
\end{equation}

Taking $x_0$  as a parameter  instead  of $s$  and substituting   (\ref{eq:2.6}) into  (\ref{eq:2.4})  we obtain  (cf. [3])
\begin{multline}																	\label{eq:2.8}
\frac{d\rho^\pm}{dx_0}=
\frac
{(f-1)\xi_0-f\xi_\rho^\pm}
{(2-f)\xi_0+(f-1)\xi_\rho^\pm}
=\frac
{(f-1)\xi_0-(f-1)\xi_0\pm \sqrt \Delta}
{(2-f)\xi_0+(f-1)\frac{-(f-1)\xi_0\pm\sqrt \Delta}{-f}}
\\
=\frac{\pm\sqrt \Delta(-f)}{-\xi_0\pm (f-1)\sqrt \Delta}
=\frac{\pm\sqrt\Delta f}{\xi_0\mp(f-1)\sqrt\Delta},
\end{multline}
where
\begin{equation}																\label{eq:2.9}
\Delta=\xi_0^2-f(\rho)\frac{\xi_\varphi^2}{\rho^2}.
\end{equation}
Analogously
\begin{equation}																\label{eq:2.10}
\frac{d\varphi^\pm}{dx_0}=\frac{-\frac{\xi_\varphi}{\rho^2}}{(2-f)\xi_0+(f-1)\xi_\rho^\pm}=
\frac{\frac{\xi_\varphi}{\rho^2}f}{-\xi_0\pm (f-1)\sqrt\Delta}=
\frac{-\frac{\xi_\varphi}{\rho^2}f}{\xi_0\mp (f-1)\sqrt\Delta}
\end{equation}
We assume that $\xi_0>0$  and  $\Delta>0$,
and we denote  by  $\rho^\pm(x_0),\varphi^\pm(x_0)$  solutions  of  (\ref{eq:2.8}),  (\ref{eq:2.10})  corresponding to the sign 
$\pm$  in   (\ref{eq:2.8})   or (\ref{eq:2.10}).

Note that  $f(\rho)=\frac{(\rho-r_+)(\rho-r_-)}{\rho^2}$.   Thus
$f(\rho)<0$  between  $\rho=r_-$  and  $\rho=r_+$   and $\Delta >0$   there.
By the  continuity  $\Delta>0$  is a small neighborhood  of  $[r_-,r_+]$.   When  $\rho\rw 0$  $\frac{f(\rho)}{\rho^2}\rw  +\infty$.
Let   $r_0$   be  such   that  $\Delta =\xi_0^2-\frac{f(r_0)}{r_0^2}\xi_\varphi^2=0$.  Thus  $\Delta>0$  for  $\rho>r_0$  and  
$\Delta <0$   for  $\rho<r_0$.

Consider  the ``minus" null-bicharacteristic
\begin{equation}																\label{eq:2.11}
\frac{d\rho^-}{dx_0}=\frac{-\sqrt \Delta  f}{\xi_0 +(f-1)\sqrt \Delta}
\end{equation}
with the initial conditions  $(\rho_0,\varphi_0,\xi_0,\xi_\varphi)$,  where  $(\rho_0,\varphi_0)$  is the point  outside  
the outer  horizon $\rho=r_+$,   i.e.   $\rho_0>r_+$.   Note that  
$$
\xi_0-\sqrt\Delta=
\frac{\xi_0^2-\big(\xi_0^2-f\frac{\xi_\varphi^2}{\rho^2}\big)}{\xi_0+\sqrt\Delta}=\frac{f\frac{\xi_\varphi^2}{\rho^2}}{\xi_0+\sqrt\Delta}.
$$
Therefore,  cancelling  $f$  we obtain
$$
\frac{d\rho^-}{dx_0}=  \frac{-f\sqrt\Delta}{f\sqrt\Delta  +\frac{f\xi_\varphi^2}{\rho^2(\xi_0+\sqrt\Delta)}}=
 \frac{-\sqrt\Delta}{\sqrt\Delta  +\frac{\xi_\varphi^2}{\rho^2(\xi_0+\sqrt\Delta)}},   \ \ \ \ \xi_0>0.
 $$
  Therefore  
   $\frac{d\rho^-}{dx_0}<0$,  i.e.  $\rho^-(x_0)$  decreases    when  $x_0$  increases.  Since  
    $\frac{d\rho^-}{dx_0}<0$  the null-bicharacterictic  $\rho^-(x_0)$  crosses  the outer  horizon  and    the inner  horizon 
    when  $x_0$  increases.
    
    Note that
    $$
    \frac{d\varphi^-}{dx_0}=-\frac
    {\frac{\xi_\varphi}{\rho^2}f}
    {\xi_0-\sqrt\Delta+f\sqrt\Delta}=
-\frac
{\frac{\xi_\varphi}{\rho^2}f}
{f\sqrt\Delta+\frac{f\frac{\xi_\varphi^2}{r^2}}{\xi_0+\sqrt\Delta}}
=
-\frac
{\frac{\xi_\varphi}{\rho^2}}
{\sqrt\Delta+\frac{\frac{\xi_\varphi^2}{\rho^2}}{\xi_0+\sqrt\Delta}}.
$$
Thus  
$\frac{d\varphi^-}{dx_0}>0$  if  $\xi_\varphi <0$  and  $\frac{d\varphi^-}{dx_0}<0$  if  $\xi_\varphi>0$.   Also  we have
\begin{equation}														\label{eq:2.12}
\frac{d\varphi^-}{d\rho^-}=\frac{\frac{\xi_\varphi}{\rho^2}}{\sqrt \Delta}.
\end{equation}
Thus 
$\frac{d\varphi^-}{d\rho^-}<0$  if  $\xi_\varphi <0$  and  $\frac{d\varphi^-}{d\rho}>0$  if  $\xi_\varphi>0$. 

Let  $r_0=r_0\big(\frac{\xi_\varphi^2}{\xi_0^2}\big)$  be  the root of   
\begin{equation}														\label{eq:2.13}
\Delta(r_0)=\xi_0^2-f(r_0)\frac{\xi_\varphi^2}{r_0^2}=0.
\end{equation}
Thus   $\Delta(\rho)=(\rho-r_0)\Delta_1(\rho^-),$            
where  $\Delta_1(\rho)>0$.   We have near  $\rho=r_0, \rho>r_0$:
$$
\frac{d\rho^-(x_0)}{dx_0}=-\sqrt{\rho^--r_0}\,\Delta_2(\rho),
$$
where  $\Delta_2(\rho)>0$.   Therefore
$$
(\rho^-(x_0)-r_0)^{\frac{1}{2}}=C(\rho)(t_0-x_0)\ \ \ \mbox{for}\  \ \ x_0<t_0.
$$
Thus  $\rho^-(x_0)$  reaches  $r_0$   when  $x_0\rw t_0,  x_0<t_0$,   i.e. 
  $\rho^-(t_0)=r_0$.   Analogously,
\begin{equation}												\label{eq:2.14}
\frac{d\varphi^-}{d\rho^-}=\frac
{\Delta_3(\rho^-)\xi_\varphi}
{(\rho^--r_0)^{\frac{1}{2}}
},
\end{equation}
where
$\rho>r_0,\Delta_3(\rho)>0$.   Hence,   assuming,  for  the  definitness,  that  $\xi_\varphi>0, $  we 
get,   for  $\rho^->r_0,\varphi^-\leq \theta_0,$
\begin{equation}												\label{eq:2.15}
(\rho^--r_0)^{\frac{1}{2}}=C_3(\rho)(\theta_0-\varphi^-),
\end{equation}
where  $r_0=\rho^-(\theta_0)$.

Note that  $\frac{d\rho^-(t_0)}{dx_0}=0$.   We shall  show  that  $(r_0,\theta_0)$  is  a turning point  of the null-geodesic
  $\gamma_0$.
For  $x_0\geq  t_0$  we  consider  the  $(+)$   solution  (cf.  (\ref{eq:2.8}))
\begin{equation}												\label{eq:2.16}
\frac{d\rho^+(x_0)}{dx_0}=\frac{f\sqrt\Delta}{\xi_0+\sqrt\Delta-f\sqrt\Delta}.
\end{equation}
Note that
$\frac{d\rho^+(t_0)}{dx_0}=0$   and       $\frac{d\rho^+(x_0)}{dx_0}>0$  for  $t_0<x_0$.

Also we have for  $x_0>t_0$
\begin{equation}												\label{eq:2.17}
\frac{d\varphi^+}{d\rho^+}=-\frac{\xi_\varphi}{\rho^2\sqrt\Delta}.
\end{equation}
Hence  $(\rho^+-r_0)^{\frac{1}{2}}=C_3(\rho)(\varphi^+-\theta_0),\ \ \ \varphi^+>\theta_0$  (cf.  (\ref{eq:2.15})).
  As above,  we are assuming that  $\xi_\varphi>0$.   Therefore
in the two-sided  neighborhood  of  $\theta_0$  we have  $\rho-r_0=C(\rho)(\varphi-\theta_0)^2$.

Thus,  $(r_0,\theta_0)$  is a turning point.

When  $x_0>t_0$  is increasing   $\rho^+(x_0)$   is also  increasing   since  $\frac{d\rho^+}{dx_0}>0$.   
Near  the inner  horizon  we have
\begin{equation}												\label{eq:2.18}
\frac{d\rho^+}{dx_0}=\frac{f\sqrt\Delta}{\xi_0+\sqrt\Delta-f\sqrt\Delta}=(\rho^+-r_-)\Delta_4(\rho^+),
\end{equation}
where  $\Delta_4(\rho^+)<0$.   Therefore,   $\frac{d\rho^+}{d x_0}\leq  -C(\rho^+-r_-)$  and integrating we  get  
$0<r_- -\rho^+(x_0)\leq  C_3  e^{-C_2 x_0}$.
Thus  $\rho^+(x_0)$  tends  to $r_-$  when  $x_0\rw+\infty$.

More precisely,  since  $f=\frac{(\rho-r_-)(\rho-r_+)}{\rho^2}$  and  $\Delta=\xi_0^2-\frac{f(\rho)\xi_\varphi^2}{\rho^2}$  we 
can  rewrite  (\ref{eq:2.18})  in the form
\begin{equation}												\label{eq:2.19}
\frac{d\rho^+}{dx_0}=\frac{(\rho^+-r_-)(r_--r_+)}{2r_-^2}+O\big((\rho^+-r_-)^2\big).
\end{equation}
Therefore
\begin{equation}												\label{eq:2.20}
r_--\rho^+=e^{-\frac{r_+-r_-}{2r_-^2}x_0+C}\big(1+O(e^{-C_1x_0})\big).
\end{equation}

We shall  summarize  the results  of this section  in the following theorem
\begin{theorem}														\label{theo:2.1}
Any ``minus"  null-geodesic  $\gamma_0$  starting above  the  outer horizon  $r=r_-$  decreases,  
 i.e.  $\frac{d\rho^-}{dx_0}<0$,
when  the time  $x_0$  increases.  It passes  the outer  and the inner horizons  $\rho=r_+$  and  $\rho=r_-$  until  
 it reaches  the turning  point  $(r_0,\theta_0),r_0<r_-$.   Then it  increases  when  the time  is increasing  and  tends  to 
the inner horizon  when  $x_0\rw +\infty$  (see  Fig.1).
\end{theorem}

{\bf Remark  2.1}
When  the  initial  point $(\rho_0,\varphi_0)$   is far  from the outer event  horizon  $\rho=r_+$  the equation   (\ref{eq:2.3})  has two roots  
$\xi_0>0$  and 
$\xi_0<0$.   We shall show that   $\rho^+(x_0)$  and $ \rho^-(x_0)$  interchange when  we change $\xi_0>0$  to $\xi_0<0$.
We  have  (cf.  (\ref{eq:2.11})):
\begin{equation}														\label{eq:2.21}
\frac{d\rho^+}{dx_0}=\frac{\sqrt \Delta f}{\xi_0-(f-1)\sqrt\Delta}=\frac{\sqrt\Delta f}{\xi_0+\sqrt\Delta -f\sqrt\Delta}.
\end{equation}
When  $\xi_0<0$  we have  $\xi_0+\sqrt\Delta=\frac{\xi_0^2-\Delta}{\xi_0-\sqrt\Delta}
=\frac{f\frac{\xi_\varphi^2}{\rho^2}}{\xi_0-\sqrt\Delta}$.
Therefore   cancelling $f$  we get
\begin{equation}														\label{eq:2.22}
\frac{d\rho^+}{dx_0}=\frac{\sqrt\Delta}{\frac{\xi_\varphi^2}{\rho^2(\xi_0-\sqrt\Delta)}-\sqrt\Delta}<0,
\end{equation}
i.e.  $\frac{d\rho^+}{d x_0}<0$.   Thus   $\rho=\rho^+(x_0)$   is  decreasing  when  $x_0$  is  increasing  and  the null-bicharachteristic  $\gamma_0$
crosses  the event  horizons  $\rho=r_+$  and   $\rho=r_-$  and   reachies  the  turning   point  $\rho=r_0$.   
 The null-bicharacteristic  $\gamma_0$ after  passing the turning point  is described by the equations  $\rho=\rho^-(x_0),\varphi=\varphi^-(x_0)$
(cf.  (\ref{eq:2.16}),  (\ref{eq:2.17})),
  and
 $$
 \frac{d\rho^-}{dx_0}=\frac{-\sqrt\Delta f}{\xi_0+(f-1)\sqrt\Delta}=\frac{-\sqrt\Delta\frac{(\rho-r_+)(\rho-r_-)}{2r^2}}{\xi_0-\sqrt\Delta +f\sqrt\Delta}>0
 $$
 for   $\rho< r_-$.  As   in  (\ref{eq:2.18})   $  \rho=\rho^-(x_0)$  tends  to  $r_-$  when  $x_0\rw +\infty$.
\\
\

\begin{tikzpicture}
[scale=2]

\draw[->] circle(0.5) (0,0) ; 

\draw[->] circle(1) (0,0) ; 

\draw[->] circle(1.5) (0,0) ; 

\draw[->] circle(0.3pt) (0,0);

\draw[->] circle(0.001pt) (-3,0);

\draw[ultra thick](-1,1.23).. controls (-1,1.15)  and (-0.8,0.6) .. (0,0.52);

\draw[ultra thick](0,0.52).. controls (0.2,0.53)  and (0.25,0.7) .. (0.3,0.9);

\draw (0.1,-1.15) node {$r=r_-$};
\draw (0.1,-1.65) node {$r=r_+$};
\draw (0.1,-0.65) node {$r=r_0$};

\end{tikzpicture}
\\
\
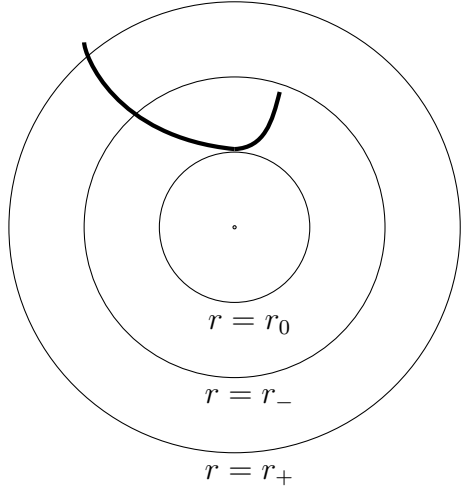
{\bf Fig. 1.}   The ``minus"  null-geodesics  crosses  the outer  and inner horizons,  makes  a turn  at some  point  
$(r_0,\theta_0),r_0<r_-,$  then  it increases  and  tends  to  the inner   horizon  when  $x_0\rw +\infty$.

\section{Geometric  optics  type  solution}
\init

The equation  $\Box_g u=0$  has  the following   form  
 in  Cartesian  coordinates  (cf.  (\ref{eq:1.3}))
\begin{equation}														\label{eq:3.1}
\frac{\partial^2 u}{\partial  x_0^2}-\sum_{k=1}^3\frac{\partial ^2 u}{\partial x_k^2}-
\Big(-\frac{\partial}{\partial x_0}
+\sum_{k=1}^3x_k\frac{\partial}{\partial  x_k}\Big)
(f-1)
\Big(-\frac{\partial}{\partial x_0}
+\sum_{k=1}^3x_k\frac{\partial}{\partial  x_k}\Big)u
=0,
\end{equation}
where   $f(r)=1-\frac{2m}{r}+\frac{e^2}{r^2}$.

Denote  by  $\Pi_0$  the plane  $x_3=0$.   The restriction  of (\ref{eq:3.1})  to the plane  $x_3=0$ 
has  the form  in polar coordinates  $(\rho,\varphi)$   (cf.  \S  2)
\begin{equation}														\label{eq:3.2}
\frac{\partial^2 v}{\partial  x_0^2}-\frac{\partial^2 v}{\partial  \rho^2}-\frac{1}{\rho^2}\frac{\partial^2 v}{\partial  \varphi^2}
-\Big(-\frac{\partial}{\partial x_0}
+\frac{\partial}{\partial  \rho}\Big)
\big(f(\rho)-1\big)
\Big(-\frac{\partial}{\partial x_0}
+\frac{\partial}{\partial  \rho}\Big)v
=0,
\end{equation}
As in \S 2  we choose  arbitrary  point  $(\rho_0,\varphi_0)$  in  the  plane  $x_3=0$  where  $\rho_0>r_+$.  Denote by  
$\gamma_0=\gamma_0(\rho_0,\varphi_0,\xi_0,\xi_\varphi)$   the null-characteristic  in the plane  $x_3=0$  starting  at  
$(\rho_0,\varphi_0,\xi_0,\xi_\rho^{(0)},\xi_\varphi)$  where  (cf.  (\ref{eq:2.7}))
\begin{equation}																\label{eq:3.3}
\xi_\rho^{(0)}= \xi_\rho^-(\rho_0,\xi_0,\xi_\varphi).
\end{equation}

Let  $V_0$  be  a neighborhood  of  $(\rho_0,\varphi_0)$  in  the plane  $x_3=0$
such  that  $|\rho'-\rho_0|<\e,|\varphi'-\varphi_0|<\e$  when  $(\rho',\varphi')\in V_0$.
Denote  by  $\gamma'(\rho',\varphi',\xi_\rho^-(\rho',\xi_0,\xi_\varphi),\xi_0,\xi_\varphi)$  the null-geodesic starting  at 
 $(\rho',\varphi',\xi_\rho^-(\rho',\xi_0,\xi_\varphi),\xi_0,\xi_\varphi)$.
  Let
\begin{equation}																\label{eq:3.4}
S_0^-(\rho',\varphi',\xi_0,\xi_\varphi)=
\int\limits_{\rho_{10}}^{\rho'}\xi_\rho^-(\rho_1,\xi_0,\xi_\varphi)d\rho_1+\varphi'\xi_\varphi,
\end{equation}
$\rho_{10}>\rho'$.

Denote  by  $\rho=\rho(x_0,\rho',\varphi'),\  \varphi=\varphi(x_0,\rho',\varphi')$  the solution  of the bicharacteristic  system  
(\ref{eq:2.8}),   (\ref{eq:2.10})  with the initial conditions  
$\rho', \varphi', \xi_\rho^-(\rho',\xi_0,\xi_\varphi),\xi_0,\xi_\varphi$.
 We have that the Jacobian 
\begin{equation}																\label{eq:3.5}
J(x_0)=\begin{vmatrix}  \frac{\partial \rho}{\partial \rho'}    &  \frac{\partial \varphi}{\partial \rho'}\\
\frac{\partial\rho}{\partial \varphi'}   & \frac{\partial\varphi}{\partial\varphi'}
\end{vmatrix}
\end{equation}
is not  zero  on  $[0,t_0-\e]$.   Therefore  the inverse map
$\rho'=\rho'(x_0,\rho,\varphi), \varphi'=\varphi'(x_0,\rho,\varphi)$ exists.   It follows from  [2],  \S 64,  that
\begin{equation}																\label{eq:3.6}
S^-(x_0,\rho,\varphi)=S_0^-(\rho'(x_0,\rho,\varphi),\varphi'(x_0,\rho,\varphi),\xi_0,\xi_\varphi)
\end{equation}
satisfies the  eikonal  equation  for  (\ref{eq:3.2})
\begin{equation}                           													\label{eq:3.7}
(2-f)(S_{x_0}^-)^2+2(f-1)S_\rho^- S_{x_0}^- -f(S_\rho^-)^2-\frac{1}{\rho^2}(S_\varphi^-)^2=0
\end{equation}
on $[0,t_0-\e]$.
Denote by  $\Pi_\alpha$   the plane passing  thrbough  the axis $0x_2$  and  having angle  $\alpha$  with the   plane  $\Pi_0$,
$\alpha\in  (-\delta,\delta)$.   
The orthogonal transformation  $O_\alpha$ 
\begin{align}																	\label{eq:3.8}
& x_1(\alpha)=x_1\cos\alpha-x_3\sin\alpha,
\\
\nonumber
& x_2(\alpha)=x_2,
\\
\nonumber
& x_3(\alpha)=x_1\sin\alpha + x_3\cos\alpha
\end{align}
maps  $\Pi_0$  onto  $\Pi_\alpha$.    Since  the Reissner-Nordstrom
metric is spherically symmetric,  the restriction  of  (\ref{eq:3.1})  to  $\Pi_\alpha$  has the
form  (\ref{eq:3.2})  in polar coordinates  in the plane  $\Pi_\alpha$  as  in the  case  $\alpha=0$.
If  $(\rho'(\alpha),\varphi'(\alpha))$  is the image of  $(\rho',\varphi')\in  V_0$   then the null-bicharacteristic  $\gamma_\alpha'$  with initial  data  $\rho'(\alpha),\varphi'(\alpha),\xi_\rho^-(\rho'(\alpha),\xi_0,\xi_\varphi),\xi_0,\xi_\varphi$   is  the image of 
$\gamma'(\rho',\varphi',\xi_\rho^-(\rho',\xi_0,\xi_\varphi),\xi_0,\xi_\varphi)$
under  the transformation  (\ref{eq:3.8}).

The  Cartesian  coordinates  of  the initial  point $(\rho',\varphi')$  in  the  plane  $x_3=0$  have  the form
$y_1^{(0)}=\rho'\cos\varphi',y_2^{(0)}=\rho'\sin\varphi',y_3^{(0)}=0$.

For  any $\alpha$  the  image  of  $(y_1^{(0)},y_2^{(0)},y_3^{(0)})$  under the map  (\ref{eq:3.8})  has the form
\begin{align}															\label{eq:3.9}
&y_1=\rho'\cos\varphi'\cos\alpha,
\\
\nonumber
&y_2=\rho'\sin\varphi,
\\
\nonumber
&y_3=\rho'\cos\varphi\sin\alpha
\end{align}
Note that the Jacobian  of the map  (\ref{eq:3.9})  is not zero  when  $(\rho',\varphi',\alpha)\subset V_0\times(-\delta,\delta)$.
Denote by  $U_0$    the image  of $V_0\times(-\delta,\delta)$  under the map  (\ref{eq:3.9}).

Let  $S_0^-(\rho',\varphi',\xi_0,\xi_\varphi,\alpha)$  be  the function   (\ref{eq:3.4})  where  $(\rho',\varphi')$
are  polar  coordinates  on the plane  $\Pi_\alpha$.   Using   the change of variable  (\ref{eq:3.9})  denote  by  
$\tilde S_0(y,\xi_0,\xi_\varphi)$  the  function  $S_0^-(\rho',\varphi',\xi_o,\xi_\varphi,\alpha)$ in the Cartesian coordinates.

  Let
\begin{equation}															\label{eq:3.10}
x=x(x_0,y)
\end{equation}
is  the  solution  of the  system  of null-bicharacteristics  with  the initial data  $x(0,y)=y$.    Since the Jacobian 
\begin{equation}               														\label{eq:3.11}
J(x_0)=\det\Big[\frac{\partial  x_j}{\partial y_k}\Big]_{j,k=1}^3
\end{equation}
is not  zero  on  $[0,t_0),  y\in U_0$,  we  have  an inverse  map
\begin{equation}															\label{eq:3.12}
y=y(x_0,x).
\end{equation}
As  in [2], \S 64,  one can  show  that
\begin{equation}															\label{eq:3.13}
S^-(x_0,x,\xi_0,\xi_\varphi)  =\tilde S_0^-(y(x_0,x),\xi_0,\xi_\varphi)
\end{equation}
is the eikonal  function  for  (\ref{eq:3.1}),  i.e. 
\begin{equation}															\label{eq:3.14}
\Big(\frac{\partial S^-}{\partial  x_0}\Big)^2-\sum_{k=1}^3\Big(\frac{\partial S^-}{\partial x_k}\Big)^2+(1-f(r))\Big(-\frac{\partial S^-}{\partial  x_0}
+\frac{\partial S^-}{\partial r}\Big)^2=0.
\end{equation}
Geometric optics solution of  (\ref{eq:3.1})  has the
following form
 form  on $ [0,t_0-\e]$  (cf.,  for  example,  [2],  \S   64)
\begin{equation}         														\label{eq:3.15}
u^-(x_0,x,k)=a_N^-(x_0,x,k)e^{ikS^-(x_0,x,\xi_0,\xi_\varphi)},
\end{equation}
where  the eikonal   $S^-(x_0,x,\xi_0,\xi_\varphi)$  satisfies  (\ref{eq:3.14}),  $k$   is  the large  parameter  and
\begin{equation}															\label{eq:3.16}
a_N^-(x_0,x,k)=a_{00}^-(x_0,x)+\frac{1}{k}a_{10}^-(x_0,x)+...+\frac{1}{k^N}a_{N0}^-(x_0,x).
\end{equation}
Note that $a_{p_0}^-(x_0,x)$  satisfy  some transport equations  (cf. [2], \S 64).  In particular,
\begin{align}																\label{eq:3.17}
&H_{0\xi_0}\Big(x,\frac{\partial  S^-}{\partial  x}\Big)\frac{\partial a_{00}^-}{\partial x_0}+\sum_{j=1}^3 H_{0\xi_j}
\Big(x,\frac{\partial  S^-}{\partial  x}\Big)\frac{\partial a_{00}^-}{\partial x_j}
\\
\nonumber
+&\Big(H_0\Big(x,\frac{\partial}{\partial x_0},\frac{\partial}{\partial x}\Big)S^- 
+f'(r)\Big(-\frac{\partial S^-}{\partial x_0}+\frac{\partial S^-}{\partial r}\Big)\Big)a_{00}^-=0.
\end{align}
Here  $H_0\Big(x,\frac{\partial}{\partial x_0},\frac{\partial}{\partial x}\Big)$  is  the principal  part  of (\ref{eq:3.1}),  $a_0^-$  
satisfies 
the initial condition
\begin{equation}															\label{eq:3.18}
a_{00}^-(0,x)=\chi(x)
\end{equation}
where  $\chi(x) \in  C_0^\infty,\ 
\mbox{supp\,}\chi\subset U_0$. 
 Note  that  $\mbox{supp\,} a_{p_0}^-(x_0,x),\ 0\leq p\leq N,$
is contained  in a neihgborhood  of  the null-geodesic $\gamma_0$.

The  geometric  optics  solution  has the form  (\ref{eq:3.15})   until  the Jacobian  (\ref{eq:3.11})  is not   zero.  
The  set  $\Sigma$   where  (\ref{eq:3.11})  is  zero   is called   the caustic set.

In  our case  the  caustic     set  $\Sigma$  is  not empty.  Its intersection  with  $\Pi_\alpha, \linebreak \alpha\in (-\delta,\delta)$, consists  of the 
circle  (\ref{eq:2.13})  
   for each  $\alpha\in(-\delta,\delta)$.    
 
Let  $x=x(x_0,y)$  be  the  equation  of null-geodesic  starting  at  $y\in \mbox{supp\,}\chi$.  Let  
$\hat a_{00}^-(x_0,y)=a_{00}^-(x_0,x(x_0,y))$.
If  we substitute  $x=x(x_0,y)$  in the equation  (\ref{eq:3.17})  we get  an ordinary  differential  equation  for
$\hat a_{00}^-(x_0,y)$   
\begin{equation}																\label{eq:3.19}
\frac{d\hat a_{00}^-(x_0,y)}{dt}=M^-(x_0,y)\hat a_{00}^-(x_0,y)
\end{equation}
where
\begin{multline}																	\label{eq:3.20} 
M^-(x_0,y)=H_{0\xi_0}^{-1}(x(x_0,y),S_{x_0}^-,S_x^-) \Big(H_0(x(x_0,y),\frac{\partial}{\partial x_0},\frac{\partial}{\partial x}\Big)S^-
\\
-f'(r)\frac{\partial S^-}{\partial  x_0}  + f'(r)\frac{\partial S^-}{\partial r}\Big).
\end{multline}
Let $(\rho',\varphi',\alpha)$  be the  pre-image of  $y=(y_1,y_2,y_3)$  under  the map  (\ref{eq:3.9}).  Then the  null-bicharacteristic  $x=x(x_0,y)$  starting  at $y$  is the  null-bicharachteristic  $\gamma'(\rho',\varphi',\xi_\rho^-(\rho',\xi_0,\xi_\varphi),\xi_0,\xi_\varphi)$   in  the plane  $\Pi_\alpha$. 
  Therefore
in $(\rho',\varphi')$  coordinates we have  $S_{x_0}^-=\xi_0,\ S_\varphi^-=\xi_\varphi$,  where  $\xi_0,\xi_\varphi$  are  constants and  
$S_\rho^-=\xi_\rho^-=\frac{-(f-1)-\sqrt\Delta}{-f}$  (cf. (\ref{eq:2.6})).
Hence substituting  in (\ref{eq:3.20})  we get  on $[0,t_0-\e]$
\begin{equation}														\label{eq:3.21}
M=\frac{(-f)\Big(-f\frac{\partial}{\partial\rho}\xi_\rho^-+f'(\rho)\xi_\rho^- -f'(\rho)\xi_0\Big)}{-\xi_0-(f-1)\sqrt\Delta}.
\end{equation}
  Note that $x_0=t_0$  is   a caustic set.  The geometric  optics  solution  $u_N^-$   is
valid  on the interval  $[0,t_0-\e]$.

When  $x_0\in  [t_0-\e,t_0+\e]$,  i.e.  in the neighborhood  of the caustic  point,  one needs to modify  the ansatz (\ref{eq:3.15})
following  Maslov theory  (cf. [6],  and also  [2],  \S 66).
  
  We will   look for modified  geometric  optics  solution   in the form
\begin{equation}														\label{eq:3.22}
u_N^{(0)}=\Big(\frac{k}{2\pi}\Big)^{\frac{1}{2}}\int\limits_{-\infty}^\infty  a_N^{(0)}(x_0,\eta,\varphi)
e^{ik\rho\eta-ikL(x_0,\eta,\varphi,\xi_0,\xi_\varphi)}d\eta,
\end{equation}
where   $L(x_0,\eta,\varphi,\xi_0,\xi_\varphi)$  satisfies   an eikonal  equation  of the form
\begin{equation}														\label{eq:3.23}
H_0\Big(-\frac{\partial L}{\partial \eta},\eta,\frac{\partial L}{\partial x_0},\frac{\partial L}{\partial \varphi}\Big)=0,
\end{equation}															
and
\begin{equation}														\label{eq:3.24}
a_N^{(0)}=\sum_{p=0}^N\frac{1}{k^p}a_{p0}^{(0)},
\end{equation}
where  $a_{p0}^{(0)}$  is satisfying   some transform  equations  (cf.  [2],  \S  66  for details).
Applying  the stationary  phase method   to the integral  (\ref{eq:3.22})  at  $x_0=t_0-\e$  
we get  a stationary  point  $\eta^{(0)}$  satisfying   the equation 
 $\rho-\frac{\partial L(t_0-\e,\eta^{(0)},\varphi,\xi_0,\xi_\varphi)}{\partial \eta}=0$.
 
 It can be shown  that  (see [2],  \S 66)
 \begin{equation}														 \label{eq:3.25}
\rho\eta^{(0)}-L(t_0-\e,\eta^{(0)},\varphi,\xi_0,\xi_\varphi)=S^-(t_0-\e,\rho,\xi_0,\xi_\varphi),
\end{equation}
where  $S^-(x_0,x,\xi_0,\xi_\rho)$  is the same as  in  (\ref{eq:3.13}).
Therefore one can adjust coefficients  $a_{p0}^{(0)}$  to have that 
\begin{equation}														\label{eq:3.26}
u_N^-\Big |_{x_0=t_0-\e}=u_N^{(0)}\Big|_{x_0=t_0-\e}
\end{equation}
modulo  lower  order  terms  in  $\frac{1}{k}$.

Analogously  applying  the stationary phase method  to  (\ref{eq:3.22})  at  $x_0=t_0+\e$  we get  modulo 
 lower  order terms  in $\frac{1}{k}$  that
\begin{equation}														\label{eq:3.27}
u_N^{(0)}\Big |_{x_0=t_0+\e}=a_N^+ e^{ikS^+(t_0+\e,\rho,\varphi,\xi_0,\xi_\varphi)},
\end{equation}
where  $S^+(x_0,x,\xi_0,\xi_\varphi)$   is the same  as  $S^-(x_0,x,\xi_0,\xi_\varphi)$  when  $\xi_\rho^-$   is  replaced  by  $\xi_\rho^+$.
Thus   for  $x_0>t_0+\e$  we  again are dealing  with the geometric  optics solutions  of the form  (\ref{eq:3.15}).

Substituting  in (\ref{eq:3.21})  $\xi_\rho^+$  instead  of  $\xi_\rho^-$  and 
using  (\ref{eq:2.20})  we  get that  $M^+$   has  the form  $M^+=O(e^{-C_1x_0})$  when  $x_0\rw +\infty$.
Therefore
solving  the ordinary  differential
 equation  (\ref{eq:3.19})  we get that  $a_{00}^+(x_0,y)$  has  a finite  limit  when  $x_0\rw +\infty$  for any fixed  $y\in U_0$.

We  shall summarize  the results  of  this section  in the following  theorem
\begin{theorem}														\label{theo:3.1}
The geometric  optics  solution  (\ref{eq:3.15})  has  the support  in  a neighborhood  of the  null-geodesoc  $\gamma_0$.  It 
has the form  (\ref{eq:3.15})  for  $x_0\leq t_0-\e$  before approaching  the caustic set.  To continue  the solution  through 
the caustic  set  one needs  to use  the  ansatz (\ref{eq:3.22}).  After  passing  the caustic  set  one  can
transform  modulo  lower  order  terms  in $\frac{1}{k}$  the ansatz   (\ref{eq:3.22})  to the    geometric  optics
ansatz  	
$u_N^+=a_N^+e^{ikS^+}$  similar  to  (\ref{eq:3.15}).  When  $x_0$  increases  $\mbox{supp\,}  a_N^+$  approaches  
the inner  horizon  $r=r_-$  and  
$\lim_{x_0\rw +\infty}a_0^+(x_0,y)$  exists  for each  $y\in U_0$.
\end{theorem}

{\bf Remark 3.1}
Let  $(\rho^{(0)},\varphi^{(0)})$    be any point  in the plane  $x_3=0$.

Assign  initial  conditions  $\rho^{(0)},\varphi^{(0)},\xi_\rho^{(0)},\xi_\varphi$  at  $x_0=0$  for the null-bicharacteristic  $\gamma^{(0)}$
where $\xi_\rho^{(0)}=\xi_\rho^-(\rho^{(0)},\xi_0,\xi_\varphi)$  is  the solution  of the quadratic equation (\ref{eq:2.6})  
with  the initial
``energy"  $\xi_0>0$  large.  As  it was  shown  in \$ 2,  $\gamma^{(0)}$  has  a turning  point  $r_0$  when
\begin{equation}																\label{eq:3.28}
\xi_0^2-f(r_0)\Big(\frac{\xi_\varphi}{r_0}\Big)^2=0.
\end{equation}
When  $\xi_0$  is large  enough  the turning  point  $r_0$  is small,   i.e.  close to the origin.

  The null-bicharacteristic  $\gamma^{(0)}$  exists  for all  $x_0>0$   and we can  also 
construct  a geometric  optics type solution $u_N$  with the support  in  a neighborhood  of $\gamma^{(0)}$  for all  $x_0>0$.  
Therefore  we can construct the solution  of initial boundary   value problem  with the support  as close to the origin  $r=0$
as we wish.  Note that the equation  (\ref{eq:3.1})  is not  hyperbolic  in a  neighborhood  of  $r=0$.   Thus not  any  initial-value   
problem  has a solution.
\\
\

{\bf Remark 3.2}

The geometric  optics solution  (\ref{eq:3.15})   is only  an approximate  solution  of the wave equation  (\ref{eq:3.1})  since
$$
\Box_g u_N=t_N e^{ikS(x_0,x,\xi_0,\xi_\varphi},
$$
where  $t_N=O\big(\frac{1}{k^N}\big)$.  
It is  possible  to find  the correction  $v_N$  such  that
$$
\Box_g v_N=-t_N e^{ikS}
$$.  
Then
$$
\Box_g(u_N+v_N)=0.
$$
To find  $v_N$  one can  use  again  geometric  optics  construction  as  in [2], \S 64  and  then  apply  the  Duhamel principle  (cf.  [4]).
    
 \section{The case  of  extremal  RN  black hole}
 \init
 
 In the case of extremal  RN  black hole,  i.e.  when  $e^2=m^2$   most of  computations  of \S2  and \S 3   remain  the same   
  with   $r_+$  and  $r_-$ replaced  by  $m$.    In particular,   formulas (\ref{eq:2.8}),  (\ref{eq:2.9}),  (\ref{eq:2.10})  hold  with  $f(\rho)=\frac{(\rho-m)^2}{\rho^2}$.
  Thus  the null-geodesic  $\gamma_0$  crosses  the event  horizon  $\{\rho=m\}$  and approaches  the turning point  $\rho=r_0$.  
  The  computations  around  the turning point are  exactly  the same  as  in \S2  (cf.  (\ref{eq:2.13})-(\ref{eq:2.17})). 
   As in  \S2  $\rho=\rho^+(x_0)$ is  increasing    after the turning point.
   
   Near the event horizon  we have  (cf.  (\ref{eq:2.18}))
   \begin{equation}															\label{eq:4.1}
   \frac{d\rho^+}{dx_0}=\frac{f\sqrt\Delta}{\xi_0+\sqrt\Delta-f\sqrt\Delta}=
   \frac{(\rho-m)^2}{\rho^2}\Delta_4(\rho),
   \end{equation}where  $\Delta_4(\rho)>0$.  Therefore
    \begin{equation}                                     											\label{eq:4.2}
    \frac{d\rho^+}{dx_0}=\frac{(\rho-m)^2}{2m^2}+O((\rho-m)^3),
    \end{equation}
   since
    \begin{equation}															\label{eq:4.3}
    \sqrt\Delta=\xi_0+O((\rho-m)^3).
    \end{equation}
    Thus 
    \begin{equation}															\label{eq:4.4}
    \frac{d\rho^+}{(\rho-m)^2+O(\rho-m)^3}=\frac{1}{2m^2}dx_0.
    \end{equation}
    Integrating,  we get
    \begin{equation}															\label{eq:4.5}
    \rho^+-m=-\frac{2m^2}{x_0}\Big(1+O\big(\frac{1}{x_0}\big)\Big).
    \end{equation}
    Compare  (\ref{eq:4.5})  and  (\ref{eq:2.20}):  In the case  of non-extremal  RN  metric  $\rho^+$  approaching   the event horizon  
    $\rho=r_-$  exponentially  when  $x_0\rw +\infty$,  but  in the  extremal  RN case  $\rho^+-m$  has only order  $O\big(\frac{1}{x_0}\big)$.
    
Now,  as in \$3,  we construct the geometric  optics  type solution  that follows  the null-bicharacteristic  $\gamma_0$ 
for the case  of 
extremal  RN  metric.  We retain  the same notations  as  in \S3.

For  $S^-(x_0,\rho,\varphi)$  (cf.  (\ref{eq:3.6}))  we have,  as in (\ref{eq:3.7}),
\begin{equation}															\label{eq:4.6}
(2-f)(S_{x_0}^-)^2  +2(f-1)S_\rho^-S_{x_0}^--f(S_\rho^-)^2-\frac{1}{\rho^2}(S_\rho^-)^2=0,
\end{equation}
where  $f=\frac{(\rho-m)^2}{\rho^2}$.   It follows  from  (\ref{eq:3.4}),  (\ref{eq:2.6})  that
\begin{multline}																\label{eq:4.7}
S_\rho^-(x_0,\rho,\varphi)=\frac{-(f-1)\xi_0-\sqrt{\xi_0^2-\frac{f\xi_\varphi^2}{\rho^2}}}{-f}=
\frac{\xi_0-\sqrt{\xi_0^2-\frac{f\xi_\varphi^2}{\rho^2}}-f\xi_0}{-f}
\\
=
-\xi_0+ \frac{-\frac{\xi_\varphi^2}{\rho^2}}{\xi_0+\sqrt{\xi_0^2-\frac{f\xi_\varphi^2}{\rho^2}}}.
\end{multline}
Thus  $S_\rho^-$  is bounded  on  $[0,t_0-\e]$  as  in \S2.   Here  $\rho^-(t_0)=r_0,\ \Delta(r_0)=0$.
Therefore  the geometric  optic  type  solution  $u^-$   on $[0,t_0-\e]$  is the same  as the corresponding  solution  in  \S3.  
On the interval  $[t_0-\e,t_0+\e]$  we also  have  the same  formulas  as in \S3.   
Consider  now the case   $x_0\geq t_0+\e$.  We have  
\begin{multline}																\label{eq:4.8}
S_\rho^+(x_0,\rho,\varphi)=\frac{-(f-1)\xi_0+\sqrt{\xi_0^2-\frac{f\xi_\varphi^2}{\rho^2}}}{-f}=
\frac{\xi_0+\sqrt{\xi_0^2-\frac{f\xi_\varphi^2}{\rho^2}}-f\xi_0}{-f}
\\
=
-\frac{2\xi_0m^2}{(\rho-m)^2}+O(\rho-m)
\end{multline}
near  $\rho=m$.   As in  (\ref{eq:3.16}),  (\ref{eq:3.17}),  (\ref{eq:3.19}),  (\ref{eq:3.20})  the principal  term  of  $u_N^+$  is  $a_{00}^+ e^{ikS^+}$
where
\begin{equation}																\label{eq:4.9}
\frac{d a_{00}^+}{dt}\equiv M^+(x_0,y)a^+_{00}(x_0,y),
\end{equation}
\begin{multline}																	\label{eq:4.10}
M^+(x_0,y)=H_{0\xi_0}^{-1}(x(x_0,y),S_{x_0}^+,S_x^+)
\\
\cdot\Big(H_0(x(x_0,y),\frac{\partial}{\partial x_0},\frac{\partial}{\partial x})S^+
-f'(r)\frac{\partial S^+}{\partial x_0}+f'(r)\frac{\partial S^+}{\partial r}\Big).
\end{multline}
In coordinates  $(\rho,\varphi,\alpha)$  we have  for  $x_0\geq  t_0+\e$
\begin{equation}																\label{eq:4.11}
H_{0\xi_0}=2(2-f)S_{x_0}^++2(f-1)S_\rho^+,
\end{equation}
where  $f=\frac{(\rho-m)^2}{\rho^2}, S_\rho^+=O\big(\frac{1}{(\rho-m)^2}\big)$  near  $\rho=m,\rho<m$.
Hence  $H_{0\xi_0}^{-1}=O((\rho-m)^2)$.   Also   $M^+=O((\rho-m)^2)$.
Taking into  account  (\ref{eq:4.5})  we get
\begin{equation}																\label{eq:412}
M^+=O\big(\frac{1}{x_0^2}\big).
\end{equation}
Therefore  $a_{00}^+$  is bounded  when $x_0\rw +\infty$  and,  as in  \S 3,   $a_{00}^+$  has  a finite  limit  when  $x_0\rw +\infty$
for any fixed $y$.
\\
\\
{\bf Remark 4.1.
The  naked singularity case.}
\\
When  $e^2>m^2$,  i.e.  $f=1-\frac{2m}{r}+\frac{e^2}{r^2}>0 $   for  all  $r>0$,  there  is  no event horizon  but the turning  point still remains  when  
$\xi_0^2-\frac{f(\rho)\xi_\varphi^2}{\rho^2}=0$.
The null-geodesic $\gamma_0$  reaches  the turning  point  at  $x_0=t_0$,  turns  and then  tends to  the infinity  when  $x_0\rw +\infty$.
The  geometric  optics  type solution  that has  the support in a neighborhood  of the null-bicharacteristic  $\gamma_0$  passes  the  caustic  set  
as  in \S3  and then  tends  to the infinity  when  $x_0\rw +\infty$,   modulo lower order  term in  $\frac{1}{k}$.

  \section{Summary}
  \init
  
  Since  the Reissner-Nordstrom  metric  is spherically  symmetric  any   plane  $\Pi_0$ passing  through
  the origin  contains null-geodesics.  We fix   such plane  and study  arbitrary  null-geodesic 
   $\gamma_0$  in this plane.  The main  feature  of
  such  null-geodesic  is that  it has  a turning point  inside  the inner  horizon.  After passing  the turning  point 
   the null-geodesic  $\gamma_0$  approaches  the inner  horizon  from the  inside   
   when  $x_0\rw +\infty$.      
  In \S  3  we use a  family  of planes  $\Pi_\alpha, |\alpha|<\delta$,  to construct  a geometric  optics  solution  $u_N$  for the wave  
  equation  (\ref{eq:1.3})  (see (\ref{eq:3.15}))   that  has the support  in a  neighborhood  of  the  null-geodesic  $\gamma_0$.
  
There is a   complication  in constructing  the geometric  optics  solution  in a neighborhood   of the  caustic set
that  requires  to use  the Maslov  theory.   We only sketch  the construction  referring  the details  to author's book  [2],  \S 66.
 
 After  passing the  caustic  set  the geometric  optics  solution  again   can  be  written  in   the  form  (\ref{eq:3.15})  and 
 it tends  to the  sphere  $r=r_-$  when  $x_0\rw +\infty$.

  In \S4  we  indicate      the changes  needed when instead  of  non-extremal RN  metrics we  consider  the case  of extremal RN  metric
  and the  case  of naked  singularity.


\begin{thebibliography}{9999}
\bibitem[5]{} A.Hamilton and P.Avelino,  The physics  of  the  relativistic  counter-streaming  instability that  drives  mass  inflation 
inside black holes,  Physics Reports  495, 1-32
\end{thebibliography}
\end{document}